# Spectral Techniques For Solving PDE Stability Model Of Vortex Rope

DIANA ALINA BISTRIAN
Department of Electrical Engineering and Industrial Informatics
Engineering Faculty of Hunedoara,"Politehnica" University of Timisoara
Str. Revolutiei Nr.5, Hunedoara, 331128
ROMANIA
diana.bistrian@fih.upt.ro

*Abstract:* In this paper spectral methods are applied to investigate the hydrodynamic instability of swirling flow with application to Francis hydraulic turbine. Spectral methods imply representing the problem solution as truncated series of smooth global functions. An $L^2$ - projection and the collocation methods are developed assessing both analytically methodology and computational techniques using symbolic and numerical conversions. Remarks concerning the efficiency and the accuracy of each method in this case are presented. The model of the trailing vortex is used to validate the numerical algorithms with existing results in the literature. All the results are compared to existing ones and they prove to agree quite well. The advantages of using this methods in flow control problems are pointed out.

## 1 Introduction

Swirling flows behaviour has long been an intensive subject of research, especially in Francis hydropower turbine design because, at partial loads, the *rope* hydrodynamic instability appears in the draft tube. The amount of computational resources required to accurately simulate the vortex rope is huge, so a complementary stability analysis is a critical requirement to predict the flow dynamics. Pozrikidis [1] offer an introductory course in fluid mechanics, covering the traditional topics in a way that unifies theory, computation, computer programming and numerical simulation. An experimental investigation of the suction side boundary layer of a large scale turbine cascade has been performed by Simoni et al. in [2]. Resiga et al. [3] carried out an experimental and theoretical investigation of the flow at the outlet of a Francis turbine runner, in order to elucidate the causes of a sudden drop in the draft tube pressure recovery coefficient at a discharge near the best efficiency operating point.

The main goal of this paper is to develop a methodology for analyzing the rope instability downstream the Francis turbine runner by means of linear stability analysis. Assessing both analytically methodology and numerical methods, the study involves new mathematical models and simulation algorithms that translate the hydrodynamical model into computer code instructions immediately following problem formulation. Classical vortex problems were chosen to validate the code with the existing results in the literature.

The hydrodynamic stability model developed in the forthcoming sections involves spectral differentiation operators derived by means of shifted orthogonal expansions of the perturbation field. The sophisticated boundary conditions corresponding to the real flow case in a Francis turbine runner motivated the use of this method, suitable for non-periodic problems with complicated boundary conditions.

Dongara [4] used the Chebyshev tau method to examine in detail a variety of eigenvalue problems arising in hydrodynamic stability studies, particularly those of Orr-Sommerfeld type. The orthogonality of Chebyshev functions was used by Bourne [5] to rewrite the differential equations as a generalized eigenvalue problem, assembling a very efficient projection based technique.

The Nobel laureate Chandrasekhar [6] presents in his study considerations of typical problems in hydrodynamic and hydromagnetic stability as a branch of experimental physics. Among the subjects treated are thermal instability of a layer of fluid heated from below, the Benard problem, stability of Couette flow, and the Kelvin-Helmholtz instability.

Many publications in the field of hydrodynamics are focused on vortex motion as one of the basic states of a flowing continuum and effects that vortex can produce. Mayer [7] and Khorrami [8] have mapped out the stability of Q-vortices, identifying both inviscid and viscous modes of instability.

The main reason for using spectral methods is their exponential accuracy. Large classes of eigenvalue problems can be solved numerically using spectral methods, where, typically, the various unknown fields are expanded upon sets of orthogonal polynomials or functions. The convergence of these methods is, in most cases, easy to assure and they are efficient, accurate and fast.

The paper is structured as follows. Some stability concepts in vortex hydrodynamics are presented in Section 2.

The hydrodynamic model of partial derivative equations (PDE) presented in Section 3 is numerically investigated developing an $L^2$-projection method and a collocation technique, both by means of operator schemes. In Section 4 considerations about these spectral methods are outlined.

The classical approaches imply a transformation of the physical domain onto the standard interval for the definition of the Chebyshev polynomials. For the approach presented in this paper, instead of using classical Chebyshev polynomials, we used shifted Chebyshev polynomials, directly defined on the physical domain of the problem $(0, r_{wall})$, preserving the orthogonality properties detailed in Section 5. The numerical algorithms was developed to work automatically for any number of expansion terms, using symbolic and numeric conversions, the implementation techniques being detailed in Section 5. In Section 6 the algorithms are validated upon the model of a trailing vortex and employed for numerical investigation of the swirl flow downstream a Francis hydraulic turbine. The main results of the paper are summarized in Section 7.

## 2 Stability Concepts In Vortex Hydrodynamics

Most of the vortex stability analyses concerned axisymmetrical vortices with axial flow [9] in order to explain the vortex breakdown phenomenon observed experimentally for the first time in pipes [10] and in hydraulic turbines [11]. Obviously, the axial symmetry hypothesis is a major simplification having the main benefit of dramatically reducing the computational cost. On the other hand, it introduces important limitations as far as the three-dimensionality and unsteadiness of the flow are concerned. Essentially, an axisymmetric flow solver provides a circumferentially averaged velocity and pressure fields that is why it is used as a basic flow for linear stability analysis.

One of the reasons to deal with a mathematical model governing the linear stability of the mechanical equilibria of a swirling flow is to investigate the evolution of the amplitudes of the velocity and the pressure perturbations fields, respectively.

In most cases, the spatially or temporal stability (classified for open flows as in [9]) under infinitesimal perturbations is reduced to the study of an algebraic eigenvalue problem which imply solving a dispersion relation connecting in fact the frequency $\omega$ and the axial wavenumber $k$ as a consequence of the condition that nontrivial eigenfunctions exist. The instability of the flow is described by the dispersion relation in the spectral space $(k,\omega)$ corresponding to the spatio-temporal evolution of the most unstable mode in the physique space $(x,t)$. Most of the investigations [8, 10, 11] concerned the values of the nondimensional parameters for which the vortex become unstable in the case of either a spatial stability or temporal stability analysis. When the complex frequency $\omega = \omega_r + i \cdot \omega_i$, $\omega_r = \text{Re}(\omega)$, $\omega_i = \text{Im}(\omega)$ is determined as a function of the real wave number $k$ a temporal stability analysis is performed. Conversely, solving the dispersion relation for the complex wave number, $k = k_r + i \cdot k_i$, $k_r = \text{Re}(k)$, $k_i = \text{Im}(k)$, when $\omega$ is given real leads to the spatial branches $k(\omega, \Psi)$ where by $\Psi$ we denoted the set of all other physical parameters involved. In both cases, the sign of the imaginary part indicates the decay or either the growth of the disturbance.

Although a spatial stability analysis implies the investigation of a nonlinear eigenvalue problem, this type of analysis directly provides the frequency ranges of the most unstable modes. More than that, the spatial stability results can be directly compared to the experimental ones since usually, in experiments, an excitation is applied to a point in the flow and then, the effect of the excitation is studied as the flow evolves downstream.

In a spatial stability analysis, useful conclusions can be drawn considering $k^2$ as a function of the real frequency $\omega$. Since $k^2 < 0$ implies an imaginary eigenvalue $k$, for the flow moving downstream the current section the only physically acceptable case is the one for which the exponential factor $e^{-|k|z}$ holds. A detailed classification of the flows with respect to the sign of the $k^2$ term was given by Benjamin [12].

Linearization of axisymmetric steady flow of incompressible and inviscid fluid using linearized Bragg-Hawthorn equation in order to analysis the stability of swirling flow downstream to the Francis turbine was used by Resiga et al [3]. In this case, the problem for the streamfunction is defined then the eigenvalue as the axial wavenumber k is obtained. Applying to the real axial and circumferential velocity profiles downstream a Francis runner at different operating points the swirling flow stability is evaluated. Following Benjamin's theory of finite transitions between frictionless cylindrical flows, an eigenvalue analysis of the linearized problem was performed. Our methods are first validated on a benchmark model of the Batchelor vortex and then applied to a real flow configuration at the outlet of a Francis turbine runner. In a preliminary analysis viscous losses can be considered negligible at the design operating point so the inviscid fluid assumption is taken. The simple stability analysis carried out in Resiga et al [3] can be recovered as a particular case here for $m=0$ and $\omega=0$.

## 3 Mathematical Model In Matrix Operator Formulation

Swirling flows models have been assessed in literature with applications to various optimization and fluid motion control problems. The hydrodynamics of rotating machines where confined vortices are developed due to the turbine rotation have been investigated in various surveys [13-16]. An experimental investigation of the suction side boundary layer of a large scale turbine cascade has been performed in [2] to study the effect of Reynolds number on the boundary layer transition process at large and moderate Reynolds numbers. The boundary element approach is assessed in [15] for the problem of the compressible fluid flow around obstacles. The system is analyzed with respect to different operating conditions, for understanding its behaviour. In [16] oscillations and rotations of a liquid droplet are simulated numerically using the level set method, and the combined effects of oscillation amplitude and rotation rate on the drop-shape oscillation is studied.

The mathematical model governing the linear spatial stability of the fluid system downstream the Francis turbine runner, corresponding to the values of tangential wavenumber $m=\pm 1$, in operator form is

$$\Lambda_1(k,F,G,H,P) = \frac{d}{dr}(rG) + rkF + mH, \quad (1)$$

$$\Lambda_2(k,F,G,H,P) = kUG - \omega G + \frac{mWG}{r} + \frac{2WH}{r} - \frac{d}{dr}P, \quad (2)$$

$$\Lambda_3(k,F,G,H,P) = kHU - \omega H + \frac{m}{r}[HW+P] + \frac{WG}{r} + G\frac{d}{dr}W, \quad (3)$$

$$\Lambda_4(k,F,G,H,P) = kFU - \omega F + \frac{FmW}{r} + G\frac{d}{dr}U + kP, \quad (4)$$

where $F$, $G$, $H$, $P$ represent the complex amplitudes of the perturbations, $k$ is the complex axial wavenumber, $m$ is the integer tangential wavenumber, $w$ represents the temporal frequency, $U$ and $W$ represent the axial and the tangential velocity, respectively, both depending only on the radial coordinate $r$.

For a given real $\omega$, the system (1)-(4) is equivalent to the complex eigenvalue problem $\Lambda_1 = \Lambda_2 = \Lambda_3 = \Lambda_4 = 0$, on the domain $(0, r_{wall})$ together with the boundary conditions in axis origin

$$H \pm G = 0, F = P = 0 \quad (5)$$

and the wall boundary conditions

$$\frac{2W_{rwall}H}{r_{wall}} - P' = 0, G = 0, \quad (6)$$

$$r_{wall}H(kU_{rwall} - \omega) \pm HW_{rwall} \pm P = 0 = 0, \quad (7)$$

$$r_{wall}F(kU_{rwall} - \omega) \pm FW_{rwall} + kr_{wall}P = 0. \quad (8)$$

associated to the bending modes investigation. The hydrodynamic eigenvalue problem (1) – (4) is written for the inviscid case in the operator form

$$k\,\Xi\mathbf{u} = \Psi\mathbf{u}, \quad \mathbf{u} = (F \ G \ H \ P)^T \quad (9)$$

where the matrix operator $\Xi$ has the nonzero elements

$$\Xi_{11} = r, \quad \Xi_{22} = U, \quad \Xi_{33} = rU, \quad \Xi_{41} = U, \quad \Xi_{44} = 1 \quad (10)$$

and the matrix operator $\Psi$ is defined as

$$\Psi_{11} = 0, \quad \Psi_{12} = -1 - d_r, \quad \Psi_{13} = m, \quad \Psi_{14} = 0, \quad (11)$$

$$\Psi_{21} = 0, \quad \Psi_{22} = \omega - \frac{m}{r}W, \quad \Psi_{23} = -\frac{2}{r}W, \quad \Psi_{24} = d_r, (12)$$

$$\Psi_{31} = 0, \quad \Psi_{32} = -W - r\frac{dW}{dr}, \quad \Psi_{33} = r\omega - mW, \quad \Psi_{34} = -m, (13)$$

$$\Psi_{41} = \omega - \frac{m}{r}W, \quad \Psi_{42} = -\frac{dU}{dr}, \quad \Psi_{43} = 0, \quad \Psi_{44} = 0, (14)$$

where $d_r$ denotes the radial derivative operator. The operator notation is widely used in the control theory community to describe and analyze systems of differential-algebraic equations. In operator formulation described in detail in [17, 18], the differential radial operator is preceded by a square, possibly singular matrix $d_r\Phi \equiv D^{(r)}\overline{\Phi}$, where $D^{(r)}$ is the spectral differentiation matrix and $\overline{\Phi}$ represents the modal collocated values of the unknown functions. The sign of the imaginary part of $k$ indicates the decay or either the growth of the disturbance. Here the flow is considered unstable when the disturbance grows, i.e. the imaginary part of $k$ is negative.

## 4 Considerations About The Spectral

# Methods

Spectral methods are one of the most used technique for the numerical investigations in hydrodynamic stability problems. Started with Orszag [19], who first used the Chebyshev spectral methods for solving hydrodynamic stability problems, many other researchers have demonstrated the applicability of this technology with high degree of accuracy: M. Khorrami, M. Malik and R. Ash [20], L. Parras and R. Fernandez-Feria [21], J. Hesthaven, S. Gottlieb and D. Gottlieb [22], Canuto et al. [23]).

The pseudospectral collocation method is associated with a grid, that is a set of nodes and that is why it is sometimes referred to as a *nodal* method. The unknown coefficients in the approximation are then obtained by requiring the residual function to be zero exactly at a set of nodes. The set of the collocation nodes is related to the set of basis functions as the nodes of the quadrature formulae which are used in the computation of the spectral coefficients from the grid values.

Instead of representing the unknown function through its values on a finite number of grid points as doing in finite difference schemes, in spectral methods the coefficients $\{f_i, g_i, h_i, p_i\}_{i=0..N}$ are used in a finite basis of known functions $\{\Phi_i\}_{i=0..N}$

$$(F, G, H, P) = \sum_{i=0}^{N}(f_i, g_i, h_i, p_i)\Phi_i \quad (15)$$

The decomposition (15) is approximate in the sense that $\{\Phi_i\}_{i=0..N}$ represent a complete basis of finite-dimensional functional space, whereas $(F, G, H, P)$ usually belongs to some other infinite-dimensional space. Moreover, the coefficients $(f_i, g_i, h_i, p_i)$ are computed with finite accuracy. Among the major advantages of using spectral methods is the rapid decay of the error, often exponential $e^{-N}$ for well-behaved functions.

## 4.1 The $L^2$ - Projection Method

Historically, this was the first method of spectral type used for nonperiodic problems.

Considering a system of partial derivative equations (PDE) in operator form
$$Lu = f, \quad (16)$$
where $L$ is the differential operator, $u$ is the vector of unknown functions, in the interval $I = (a,b)$, coupled with the boundary conditions
$$u(a) = \lambda_1, \quad u(b) = \lambda_2. \quad (17)$$
the PDE system is required to be satisfied at each point in its domain. We introduce a finite basis $\{\Phi_i\}_{i=0..N}$ of orthogonal polynomials with respect to a weight function $w$ in the Hilbert space $L_w^2$, which satisfy $\deg \Phi_i = i$ and $(\Phi_i, \Phi_j)_w = c_i \delta_{ij}$ $(i,j = 0,1,...)$ for suitable constants $c_i > 0$. Examples are the Chebyshev system $\{T_i, i = 0,1,...\}$, for which $w(x) = (1-x^2)^{-1/2}$, the Legendre system $\{L_i, i = 0,1,...\}$, for which $w(x) = 1$, or, more generally, any Jacobi system $\{P_i^{(\lambda,\mu)}, i = 0,1,...\}$, for which $w(x) = (1-x)^\lambda (1+x)^\mu$, $\lambda, \mu > -1$.

The discrete solution is therefore represented as
$$u^N(x) = \sum_{i=0}^{N}\hat{u}_i \Phi_i(x) \quad (18)$$
where the unknowns are the expansion coefficients of $u^N$ along the chosen basis, computed as
$$\hat{u}_i = \frac{(u^N, \Phi_i)_w}{(\Phi_i, \Phi_i)_w} \quad (19)$$

The boundary conditions (17) impose two linear combinations upon the coefficients of $u^N$, namely
$$\sum_{i=0}^{N}\hat{u}_i \Phi_i(a) = \lambda_1, \quad \sum_{i=0}^{N}\hat{u}_i \Phi_i(b) = \lambda_2. \quad (20)$$

The residual $r(u^N) = f - Lu^N$ is required to be orthogonal to all polynomials of degree up to $N-2$, meaning that
$$(Lu^N, \Phi_j)_w = (f, \Phi_j)_w, \quad 0 \le j \le N-2. \quad (21)$$

At the algebraic level, this method produces a linear system of the form
$$M\overline{u} = \overline{f} \quad (22)$$
where $\overline{u} = (\hat{u}_0, ..., \hat{u}_N)$ is the vector collecting the unknowns that represent $u^N$, $\overline{f} = (\hat{f}_0, ..., \hat{f}_{N-2}, \lambda_1, \lambda_2)$ is a known vector depending on the data $f$ and the valued on the boundary, and $M$ is the matrix corresponding to the equations defined by the method.

## 4.2 The Collocation Method

If some of the coefficients of the equation are variable, the projection method is much less efficient and the collocation method is an efficient alternative.

Consider again a system of partial derivative equations (PDE) in operator form
$$Lu = f, \quad (23)$$
where $L$ is the differential operator, $u$ is the vector of unknown functions, in the interval $I = (a,b)$, coupled with the boundary conditions

$$u(a) = \lambda_1, \quad u(b) = \lambda_2. \quad (24)$$

The collocation method is associated with a grid of clustered nodes $x_j$ and weights $w_j$ ($j = 0,...,N$). The collocation nodes must cluster near the boundaries to diminish the negative effects of the Runge phenomenon [24]. Another aspect is that the convergence of the interpolation function on the clustered grid towards unknown solution is extremely fast.

We recall that the nodes $x_0$ and $x_N$ coincide with the endpoints of the interval $[a,b]$, and that the quadrature formula is exact for all polynomials of degree $\leq 2N-1$, i. e.,

$$\sum_{j=0}^{N} v(x_j) w_j = \int_a^b v(x) w(x) dx, \quad (25)$$

for all $v$ from the space of test functions.

Let $\{\Phi_\ell\}_{\ell=0..N}$ a finite basis of polynomials relative to the given set of nodes, not necessary being orthogonal. If we choose a basis of non-orthogonal polynomials we refer to it as a *nodal* basis (Lagrange polynomials for example).

An example of nodal basis is given by Lagrange's formula

$$\Phi_\ell(x) = \prod_{\substack{j \neq \ell \\ 0 \leq j, \ell \leq N}} \frac{(x-x_j)}{(x_\ell - x_j)} \quad (26)$$

For numerical stability reasons, often Lagrangian polynomials are reformulated in barycentric form as

$$\Phi_\ell(x) = \frac{\lambda_\ell}{x-x_\ell} \left( \sum_{k=0}^{N} \frac{\lambda_k}{x-x_k} \right)^{-1}, \lambda_\ell = \left( \prod_{k \neq \ell} \frac{1}{(x_\ell - x_k)} \right)^{-1} \quad (27)$$

In nodal approach, each function of the nodal basis is responsible for reproducing the value of the polynomial at one particular node in the interval.

A different approach is obtained by taking as basis functions simple linear combinations of orthogonal polynomials. These are called bases of *modal* type, i. e., such that each basis function provides one particular pattern of oscillation of lower and higher frequency.

Examples of simple modal bases are the following

$$\Phi_\ell(x) = T_\ell(x) - T_{\ell+2}(x), \quad \ell = 0..N, \quad (28)$$

$$\Psi_\ell(x) = T_\ell(x) - 2\frac{(\ell+2)}{(\ell+3)} T_{\ell+2}(x) + \frac{(\ell+1)}{(\ell+3)} T_{\ell+4}(x), \quad (29)$$

where $T_\ell(x)$ are the Chebyshev polynomials, or the modal basis functions

$$\Theta_\ell(x) = \sqrt{\frac{2\ell+3}{2}} \left( \frac{L_{\ell+3} - L_{\ell+1}}{(2\ell+3)(2\ell+5)} - \frac{L_{\ell+1} - L_{\ell-1}}{(2\ell+1)(2\ell-1)} \right), \quad (30)$$

used by Melenk, Kirchner and Schwab that utilized a Legendre-Galerkin approximation in [25].

Boundary-adapted bases of *modal* type are also useful in numerical approximations of hydrodynamic stability problems. From their construction they contains two functions that are nonzero at precisely one endpoint of the interval, which are called *vertex basis functions* and $N-1$ functions that vanish at both endpoints, which are called *bubble functions or internal basis functions* [23]. An example of boundary adapted modal basis is as follows

$$\Gamma_0(x) = \frac{1}{2}(\eta_0(x) - \eta_1(x)) = \frac{1-x}{2},$$

$$\Gamma_1(x) = \frac{1}{2}(\eta_0(x) + \eta_1(x)) = \frac{1+x}{2}, \quad (31)$$

$$\Gamma_\ell(x) = \begin{cases} \eta_0(x) - \eta_\ell(x), & \ell \text{ even} \geq 2 \\ \eta_1(x) - \eta_\ell(x), & \ell \text{ odd} \geq 3 \end{cases}, \quad 2 \leq \ell \leq N$$

where $\eta_\ell(x)$ denotes either $T_\ell(x)$ or $L_\ell(x)$.

A comparison of the behaviour of the members of the three bases mentioned in (28), (29) and (30) is given in Figure 1, for $N = 3$.

The efficiency of the collocation based algorithms was exposed in [26], for solving the Hartree-Fock equations of the self-consistent field in large atomic and molecular systems.

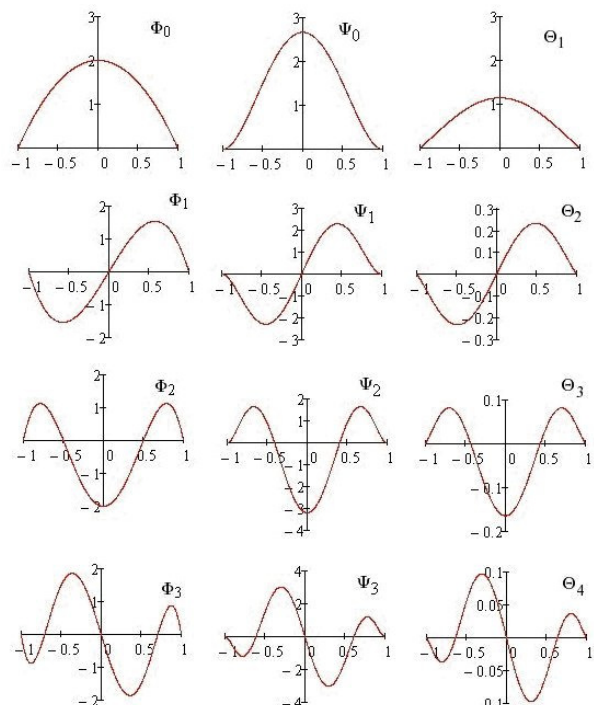

Fig.1 The modal orthogonal basis $\{\Phi_\ell\}$ given by (28) (*left*), the modal basis $\{\Psi_\ell\}$ given by (29) (*center*) and the modal basis $\{\Theta_\ell\}$ given by (30) (*right*), for $N = 3$.

# 5 Numerical Techniques For Bending Modes Investigation

## 5.1 Orthogonal Basis Functions

The numerical methods considered here are the Chebyshev pseudospectral methods. Finite element techniques reconstruct functions from a superposition of piecewise polynomial functions on subsets of triangulations of a domain or its boundary. In contrast to this, the pseudospectral techniques surveyed here will avoid triangulations and meshing, but the unknown functions are reconstructed by the superposition of simple functions.

There are two possible approaches of the mathematical model at this point. The first one imply a transformation of the physical domain onto the standard interval of the definition of the Chebyshev polynomials [27] and in the second one, instead of using classical Chebyshev polynomials, we used shifted Chebyshev polynomials $T_k^*$, directly defined on the physical interval of the problem. This second approach was our choice motivated by the form of the singular coefficients in the equations defining the eigenvalue problem and also on our previous investigation choice [28].

The shifted Chebyshev polynomials of the first kind $T_n^*(r)$ of degree $n-1$ in r on $[0, r_{wall}]$ are given by

$$T_n^*(r) = T_n\left(2rr_{wall}^{-1} - 1\right) \quad (32)$$

The shifted Chebyshev class is orthogonal in the Hilbert space $L_w^2(0, r_{wall})$, weighted by $w(r) = \left(1 - \left(2rr_{wall}^{-1} - 1\right)\right)^{-1/2}$ and have the orthogonality properties

$$\left(T_n^*, T_m^*\right)_w = 0, \quad n \neq m, n, m = 1..N, \quad (33)$$

$$\left(T_n^*, T_n^*\right)_w = r_{wall}\frac{\pi}{\lambda}, \quad \lambda = \begin{cases} 2 & \text{if } n = 1 \\ 4 & \text{if } n = 2..N \end{cases}, \quad (34)$$

with respect to the inner product $(u,v)_w = \int_0^{r_{wall}} uvw\, dr$.

For the implementation procedure, we define the derivatives by the following formulae

$$T_1^{*'} = 0$$

$$T_2^{*'} = \frac{2}{r_{wall}} T_1^*$$

$$T_3^{*'} = \frac{8}{r_{wall}} T_2^*$$

$$T_4^{*'} = \frac{6}{r_{wall}}\left[2T_3^* + T_1^*\right]$$

$$T_5^{*'} = \frac{8}{r_{wall}}\left[2T_4^* + 2T_2^*\right]$$

$$T_6^{*'} = \frac{10}{r_{wall}}\left[2T_5^* + 2T_3^* + T_1^*\right]$$

$$T_7^{*'} = \frac{12}{r_{wall}}\left[2T_6^* + 2T_4^* + 2T_2^*\right]$$

$$T_8^{*'} = \frac{14}{r_{wall}}\left[2T_7^* + 2T_5^* + 2T_3^* + T_1^*\right]$$

$$T_9^{*'} = \frac{16}{r_{wall}}\left[2T_8^* + 2T_6^* + 2T_4^* + 2T_2^*\right]$$

$$T_{10}^{*'} = \frac{18}{r_{wall}}\left[2T_9^* + 2T_7^* + 2T_5^* + 2T_3^* + T_1^*\right].$$

The four unknown components of the perturbation field are written as truncated series of orthonormal shifted Chebyshev polynomials $T_k^*$

$$(F, G, H, P) = \sum_{k=1}^{N} (f_k, g_k, h_k, p_k) \cdot T_k^* \quad (35)$$

Since, in order to discretize our hydrodynamic stability problem, a much more convenient choice is the range $[0, r_{wall}]$ than the standard definition interval of classical Chebyshev polynomials $[-1, 1]$, the independent variable $\xi \in [-1, 1]$ is maped to the variable $r \in [0, r_{wall}]$ by the linear transformation

$$r = r_{wall}(\xi + 1)2^{-1} \quad (36)$$

Consider the one dimensional domain $0 \leq r \leq r_{wall}$, where $r_{wall}$ means the radial distance to the wall. The domain of interest is represented by the Chebyshev-Gauss points in radial direction

$$\{r_k\}_{k=1}^{N} = \frac{r_{wall}}{2}\left\{1 + \cos\left(\frac{\pi(i + N - 1)}{N - 1}\right)_{i=0}^{N-1}\right\} \quad (37)$$

clustered near the boundaries.

For the case $m = \pm 1$ that we investigate here, the boundary conditions read

$$\frac{2W_{rwall}}{r_{wall}}\sum_1^N h_k - p_2\frac{2}{r_{wall}} - \sum_{3, k\,odd}^N p_k \frac{2(k-1)}{r_{wall}}\left[\sum_{r=k-1, k\,even}^2 2\right] -$$

$$- \sum_{4, k\,even}^N p_k \frac{2(k-1)}{r_{wall}}\left[\sum_{r=k-1, k\,odd}^2 2 + 1\right] = 0, \quad (38)$$

$$kU_{rwall}r_{wall}\sum_1^N h_k + (\pm W_{rwall} - \omega r_{wall})\sum_1^N h_k \pm \sum_1^N p_k = 0, \quad (39)$$

$$k\left(U_{rwall}r_{wall}\sum_1^N f_k + r_{wall}\sum_1^N p_k\right) + (\pm W_{rwall} - \omega r_{wall})\sum_1^N f_k = 0, (40)$$

$$\sum_1^N (-1)^{k+1} g_k \pm \sum_1^N (-1)^{k+1} h_k = 0, \quad \sum_1^N g_k = 0, \quad (41)$$

$$\sum_1^N (-1)^{k+1} f_k = \sum_1^N (-1)^{k+1} p_k = 0. \quad (42)$$

## 5.2 $L^2$-Projection Algorithm

The $L^2$-projection method, also known as the Chebyshev tau method, have been the attention of much study and has been successfully applied to many hydrodynamic stability problems. This represents an efficient numerical technique to solve eigenvalue problems with sophisticated boundary conditions by translate it into a linear system of equations. D. Bourne in [5] is examined the Chebyshev tau method using the orthogonality of Chebyshev functions to rewrite the differential equations as a generalized eigenvalue problem. This problem is addressed here, in application to the Benard convection problem, and to the Orr-Sommerfeld equation which describes parallel flow. J.J. Dongarra, B. Straughan and D.W. Walker in [4] examined in detail the Chebyshev tau method for a variety of eigenvalue problems arising in hydrodynamic stability studies, particularly those of Orr-Sommerfeld type. Physical problems explored in this study are those of Poiseuille flow, Couette flow, pressure gradient driven circular pipe flow, and Couette and Poiseuille problems for two viscous immiscible fluids.

Following [5] the difficult eigenvalue problem (1)-(8) is transformed into a system of linear equations describing the hydrodynamic context for the cases $m = \pm 1$. The difference between the classical tau method and the modified $L^2$-projection proposed here is given by the selected spaces involved in the approximation process. An appropriate solution is sought in the truncated Chebyshev series form (35).

The projection method is an algorithm implying in the first step expanding the residual function as a series of shifted Chebyshev polynomials. We obtain a set of $4(N-2)$ linear equations. The eight remaining equations are provided by the boundary conditions applied as side constraints.

Introducing the notations
$$I^U_{ijkld} = \left(r^k(U^l)^{(d)}T^*_i, T^*_j\right)_w, I^W_{ijkld} = \left(r^k(W^l)^{(d)}T^*_i, T^*_j\right)_w, (43)$$
with $d$ the derivation order, the first truncated $4(N-2)$ equations of the hydrodynamic model read

$$k\sum_{j=1}^N \left(f_j I^U_{ij100}\right) + g_i r_{wall} c + g_2 \frac{2}{r_{wall}} I^U_{i1100} +$$
$$+ \sum_{j=3, j\,odd}^N g_j \frac{2(j-1)}{r_{wall}} \left[\sum_{r=j-1, j\,even}^2 2 I^U_{ir100}\right] +$$
$$+ \sum_{j=4, j\,even}^N g_k \frac{2(j-1)}{r_{wall}} \left[\sum_{r=j-1, j\,odd}^2 \left(2I^U_{ir100}\right) + I^U_{i1100}\right] + m h_i r_{wall} c = 0, (44)$$

$$k\sum_{j=1}^N g_j I^U_{ij010} - \omega g_i r_{wall} c + m\sum_{j=1}^N g_j I^W_{ij-110} +$$
$$+ 2\sum_{j=1}^N h_j I^W_{ij-110} - p_2 \frac{2}{r_{wall}} A_{11} - \sum_{j=3, j\,odd}^N p_j \frac{2(j-1)}{r_{wall}}\left[\sum_{r=j-1, j\,even}^2 2A_{ir}\right] -$$
$$- \sum_{j=4, j\,even}^N p_j \frac{2(j-1)}{r_{wall}}\left[\sum_{r=j-1, j\,odd}^2 (2A_{ir}) + A_{i1}\right] = 0, (45)$$

$$k\sum_{j=1}^N h_j I^U_{ij110} - \omega \sum_{j=1}^N h_j I^U_{ij100} + m\sum_{j=1}^N h_j I^W_{ij010} + m p_i r_{wall} c +$$
$$+ \sum_{j=1}^N g_j \left(I^W_{ij010} + I^W_{ij111}\right) = 0, \quad (46)$$

$$k\left(\sum_{j=1}^N f_j I^U_{ij010} + p_i r_{wall} c\right) - \omega f_i r_{wall} c +$$
$$+ m\sum_{j=1}^N f_j I^W_{ij-110} + \sum_{j=1}^N g_j I^U_{ij011} = 0, \quad (47)$$

where the number $c$ is defined as
$c = \begin{cases} \pi/2, & j=1 \\ \pi/4, & j=2...N-2 \end{cases}$ and $A \in M_{N-2}$ is square $(N-2)\times(N-2)$ matrix, with $A_{11} = r_{wall}\pi/2$, $A_{jj} = r_{wall}\pi/4$, $j=2...N-2$, $A_{ij} = 0$, $i \neq j$. Similarly, we translate the boundary conditions in linear equations that complete the system. For the case $m = \pm 1$ that we investigate here, the boundary conditions provide the seven equations set (38) - (42).

For this case, the additional relation in obtained from the second equation of the mathematical model (2), taking the inner product $\left(\Lambda_2, T^*_{N-1}\right)_w$.

The eigenvalue problem is written as a system of $4N$ equations with the matrix formulation $kM_k \bar{s} = M\bar{s}$, $\bar{s} = \left(\bar{f}, \bar{g}, \bar{h}, \bar{p}\right)$, with $\bar{*} = (*_1, ..., *_N)$, $* \equiv f, g, h, p$. The method has the obvious advantage that the highest degree of the Chebyshev polynomials multiplying the residual in the method inner-product is only $N-2$.

The recurrence relation for $T^*_n$ has the form
$$T^*_n(r) = 2\left(2rr^{-1}_{wall} - 1\right)\cdot T^*_{n-1}(r) - T^*_{n-2}(r), n=3,4,... (48)$$
in which the initial conditions are $T^*_1(r) = 1, T^*_2(r) = \frac{2r}{r_{wall}} - 1$. The use of a recurrence relation significantly increases the elapsed time to generate the shifted Chebyshev polynomials. To improve the performance of the numerical algorithm, we introduce in our code the equivalent polynomial relation

$$T^*_n(r) = \frac{1}{2}\left[\left(\tilde{r} + \sqrt{(\tilde{r})^2 - 1}\right)^{n-1} + \left(\tilde{r} - \sqrt{(\tilde{r})^2 - 1}\right)^{n-1}\right], \tilde{r} = \frac{2r}{r_{wall}} - 1 \quad (49)$$

to automatically generate the shifted Chebyshev polynomials $\{T_n^*(\xi)\}_{n\geq 1}$ on $[0, r_{wall}]$.

The numerical algorithm was developed to work automatically for any number of expansion terms, using the routines of a high level language such Matlab.

### 5.3 Spectral Collocation Technique

The collocation method that we present in this section has the peculiar feature that can approximate the perturbation field for all types of boundary conditions, especially when the boundary limits are described by sophisticated expressions.

We assume for this approach the model and the boundary conditions described in Section 3.

This collocation method applied to stability investigation of swirling flows was applied successfully and validated in our previous study [28].

Following standard procedures, the Chebyshev spectral collocation method can be described as follows. An approximation based on Chebyshev polynomials to the unknown functions is first introduced. The set of collocation equations is then generated. The equation system consists of two parts. The first part is formed by making the associated residual equal to zero at the collocation points, while the second part is obtained by forcing the boundary conditions to be satisfied at the boundary collocation points.

In the radial direction, the values of relevant derivatives with respect to $r$ at the grid points are computed by the differentiation matrix operator $D^{(r)}$, that was derived in our previous study [18], as

$$d_r \bar{s} \equiv D^{(r)} \bar{s}, \quad \bar{s} = (f_1,...,f_N, g_1,...,g_N, h_1,...,h_N, p_1,...,p_N)^T \quad (50)$$

Let us denote by $[r] = diag(r_i)$, $\left[\frac{1}{r}\right] = diag(1/r_i)$, $[U] = diag(U(r_i))$, $[W] = diag(W(r_i))$, $2 \leq i \leq N-1$, $[\eta] = (\eta_{ij})_{\substack{2\leq i\leq N-1\\1\leq j\leq N}}$, $\eta_{ij} = T_j^*(r_i)$, $[D] = (D^{(r)}_{ij})_{\substack{2\leq i\leq N-1\\1\leq j\leq N}}$, $\bar{f} = (f_1,...f_N)$, $\bar{g} = (g_1,...g_N)$, $\bar{h} = (h_1,...h_N)$, $\bar{p} = (p_1,...p_N)$.

The eigenvalue problem (9) governing the inviscid stability analysis is written in the computational form

$$k[r][\eta]\bar{f} + ([D]+[\eta])\bar{g} + m[\eta]\bar{h} = 0, \quad (51)$$

$$\left\{k[U][\eta] - \omega[\eta] + m\left[\frac{W}{r}\right][\eta]\right\}\bar{g} + 2\left[\frac{W}{r}\right][\eta]\bar{h} - [D]\bar{p} = 0, (52)$$

$$\{[W][\eta] + [r][W'][\eta]\}\bar{g} + m[\eta]\bar{p} +$$
$$+ \{-\omega[r][\eta] + m[W][\eta] + k[r][U][\eta]\}\bar{h} = 0, \quad (53)$$

$$\left\{k[U][\eta] - \omega[\eta] + m\left[\frac{W}{r}\right][\eta]\right\}\bar{f} + [U'][\eta]\bar{g} + k[\eta]\bar{p} = 0. (54)$$

Solving the resulting eigenvalue problem with nonconstant coefficients imply imposing that equations (51) – (54) to be satisfied at the $(N-2)$ interior points $(r_i), i = 2..N-1$. The system of $4N$ equations is completed with the boundary relations (38) - (42), respectively.

## 6 Numerical Results Of The Vortex Rope Stability Investigation

### 6.1 Code Validation And Error Analysis

The basic flow under consideration for the validation of the proposed method is the Batchelor vortex case or the q-vortex [29], that trails on the tip of each delta wing of the airplanes. The properties of the Batchelor vortex were pointed out in Olendradru et al. [29] using a shooting method. In order to compare our results with the ones from [29] numerical evaluations of the axial wavenumber $k$ were obtained for various sets of parameters associated with the investigated modes. In Table 1 these values are presented in comparison with the ones from reference [29]. The numerical results obtained employing the presented methods are in agreement with the results presented in reference [29].

Table 1. Comparative results of the most amplified spatial wave of the Batchelor-vortex: eigenvalue with largest imaginary part $k_{cr} = (k_r, k_i)$.

| $m = 1$ | $m = -1$ |
|---|---|
| $L^2$ Projection method ||
| $k_{cr} = (-1.1606, 0)$ | $k_{cr} = (0.9046, -0.783)$ |
| Collocation method ||
| $k_{cr} = (0.5611, 0)$ | $k_{cr} = (0.76146, -0.33722)$ |
| Shooting method [29] ||
| $k_{cr} = (0.6, 0)$ | $k_{cr} = (0.761, -0.336)$ |

Although the collocation algorithm is a very efficient technique, the inclusion of the boundary conditions as equations in the system of the generalized eigenvalue problem (9) have been observed to be one cause of spurious eigenvalues. The spurious eigenvalues, which are not always easy to identify, may lead one to a false conclusion regarding the stability of the fluid system, thus the elimination of them is of great importance. These are values returned by the algorithm which do not satisfy the eigenvalue problem. The spurious eigenvalues problems have been the attention of much study recently. Gardner et al. [30] describe the tau methods to avoid spurious eigenvalues and in Dongara [4] the occurrence of the spurious eigenvalues is assessed in application to the Benard convection problem.

We implement in our numerical procedure a code sequence that identifies if an eigenvalue of the spectra is spurious or not. First the algorithm provides the entire spectra, then calculates the residual vector of the eigenvalue problem (9) for each eigenvalue of the spectra. A true value of $k$ must satisfy the eigenvalue problem. We evaluate the $L^2$ norm of the vector with respect to a given tolerance $\varepsilon$. If the condition

$$\left(\sum abs(k\Xi\mathbf{u}-\Psi\mathbf{u})^2\right)^{1/2} > \varepsilon, \quad \mathbf{u} = (F \quad G \quad H \quad P)^T \quad (55)$$

holds, the eigenvalue $k$ is declared spurious and discarded from the spectra.

## 6.2 Algorithm Accuracy And Numerical Results

We designed the stability algorithm in two stages. First, the algorithm solves the eigenvalue problem (9) and finds the critical eigenvalue with largest negative imaginary part, that corresponds to the most unstable perturbation. In the second stage, varying the frequency omega at different numbers of collocation parameters $N$, we retain the maximum growth rate and the corresponding frequency, denoted as the critical frequency. The question here is how to find the optimum value of the spectral parameter $N$ that defines the number of Chebyshev collocation nodes?

Let us define the eigenvalue problem (9) in operator formulation

$$\left(kL^{[k]} + \omega L^{[\omega]} + L\right)\mathbf{u} = 0, \quad \mathbf{u} = (F \quad G \quad H \quad P)^T \quad (56)$$

where

$$L^{[k]} = \begin{pmatrix} r & 0 & 0 & 0 \\ 0 & U & 0 & 0 \\ 0 & 0 & rU & 0 \\ U & 0 & 0 & 1 \end{pmatrix}, \quad L^{[\omega]} = \begin{pmatrix} 0 & 0 & 0 & 0 \\ 0 & -1 & 0 & 0 \\ 0 & 0 & -r & 0 \\ -1 & 0 & 0 & 0 \end{pmatrix},$$

$$L = \begin{pmatrix} 0 & 1+d_r & -m & 0 \\ 0 & mW/r & 2W/r & -d_r \\ 0 & W+rdW/dr & mW & m \\ mW/r & dU/dr & 0 & 0 \end{pmatrix} \quad (57)$$

Let us denote by

$$\Lambda_\omega = \left\{k = (k_r, k_i) \mid \left(kL^{[k]} + \omega L^{[\omega]} + L\right)\mathbf{u} = 0\right\} \quad (58)$$

the spectra of the eigenvalue problem (56) computed for a given frequency, and

$$gr_\omega = \min(imag\,\Lambda_\omega) \quad (59)$$

the growth rate of the most unstable perturbation at a given frequency.

Let us define the set

$$\chi_N(\omega) = \left\{-gr_\omega \mid \omega \in [0, 0.4]\right\} \quad (60)$$

and the pair

$$(gr_{\max}, \omega_{cr}) = \left\{gr_{\max} = \max(\chi_N(\omega)), \omega_{cr} = \chi_N^{-1}(gr_{\max})\right\} (61)$$

The set $\chi_N(\omega)$ and the pair (61) are computed for each mode number investigated $m = \{-1, 1\}$ for an optimum collocation number of nodes $N$, employing both projection algorithm and collocation method and are presented in Table 2.

Table 2. The critical frequency and the maximum growth rates $(\omega_{cr}, gr_{\max})$ obtained for the investigated modes, employing projection method and collocation.

| Mode number | $m = -1$ | $m = 1$ |
|---|---|---|
| $L^2$ Projection | $(0.35, 5.20004)$ | $(0.2, 0.12340)$ |
| Collocation | $(0.3, 6.93195)$ | $(0.3, 0.54437)$ |

Following Tadmor [31], when the Chebyshev pseudospectral methods are used, the error committed is expected to decay to zero at an exponential rate. For this reason, we run the algorithm for the values of collocation parameter $N$ along an interval sufficiently large to reach the convergence, for example, between 5 and 60. The output $\omega_{cr}$ that is returned over the optimum value of the collocation parameter $N_{cr}^{\min}$ is expected to have the greatest number of occurences.

The convergence behaviour of the algorithms with respect to the number of expansion terms is shown in Table 3 and Table 4.

Clearly the numerical computation costs were less expansive in the projection method approach since the number of terms in the approximations was significantly reduced. In fact, in comparison with the collocation method this number was more than two times reduced. As a result, with a reduced by far computational time, we can obtain accurate results in an acceptable agreement with existing ones.

Let us define

$$\overline{e(N)} = \left(\sum abs\left(\left(kL^{[k]} + \omega_{cr}L^{[\omega]} + L\right)\mathbf{u}\right)^2\right)^{1/2}, \mathbf{u} = (F \quad G \quad H \quad P)^T (62)$$

the residual vector of the eigenvalue problem (56) computed for a frequency equal to the critical frequency $\omega_{cr}$ and for a spectral parameter $N \in \left[N_{cr}^{\min}, N_{cr}^{\max}\right]$ along the optimum interval where the algorithm convergence is achived. Let be

$$E_N = \left\{\max\left(\overline{e(N)}\right) \middle| N \in \left[N_{cr}^{\min}, N_{cr}^{\max}\right]\right\} \quad (63)$$

the set of maximum values of the residual along the optimum interval of collocation.

Table 3. The convergence behaviour of the projection algorithm for the investigated mode numbers.

| $m = -1$ | | $m = 1$ | |
|---|---|---|---|
| N | Frequency $\omega_{cr}$ | N | Frequency $\omega_{cr}$ |
| 3 | 0.05 | 3 | 0.1 |
| 4 | 0.4 | 4 | 0.1 |
| 5 | 0.25 | 5 | 0.15 |
| 6 | 0.35 | 6 | 0.25 |
| 7 | 0.35 | 7 | 0.2 |
| 8 | 0.35 | 8 | 0.2 |
| 9 | 0.35 | 9 | 0.2 |
| 10 | 0.35 | 10 | 0.2 |

Table 4. Convergence of the collocation algorithm for the investigated mode numbers.

| $m = -1$ | | $m = 1$ | |
|---|---|---|---|
| N | Frequency $\omega_{cr}$ | N | Frequency $\omega_{cr}$ |
| 5 | 0.25 | 6 | 0.1 |
| 6 | 0.1 | 8 | 0.2 |
| 8 | 0.1 | 10 | 0.35 |
| 10 | 0.25 | 15 | 0.25 |
| 12 | 0.3 | 20 | 0.2 |
| 16 | 0.3 | 27 | 0.35 |
| 17 | 0.3 | 33 | 0.3 |
| 19 | 0.3 | 36 | 0.3 |
| 29 | 0.3 | 39 | 0.3 |
| 31 | 0.3 | 42 | 0.3 |
| 37 | 0.3 | 55 | 0.3 |
| 40 | 0.3 | 57 | 0.3 |
| 46 | 0.3 | 60 | 0.3 |
| 61 | 0.3 | 64 | 0.3 |

For each mode number we present the set values (63) on a logarithmical representation in Figures 2-3. This proves that the convergence is reached for all mode numbers considered here. One may notice that the optimum node number varies function of mode number. Using the logarithmic representation is possible to observe the fact that, along the optimum interval of collocation, the error is not decreasing with the number of nodes as expected, having a rather an exponential increase. This leads to the idea that there is an optimum number of nodes located into the optimum collocation interval, which is not necessary be a large value.

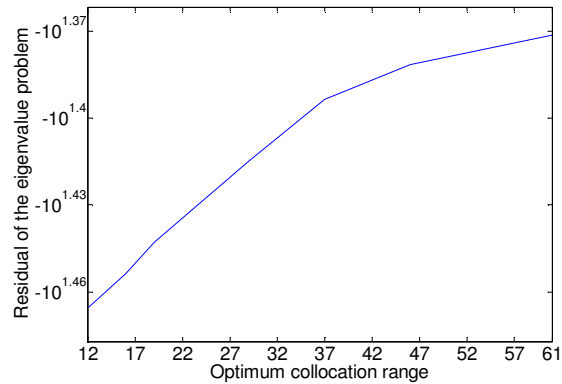

Fig.2 Residual along the optimum range for mode $m = -1$.

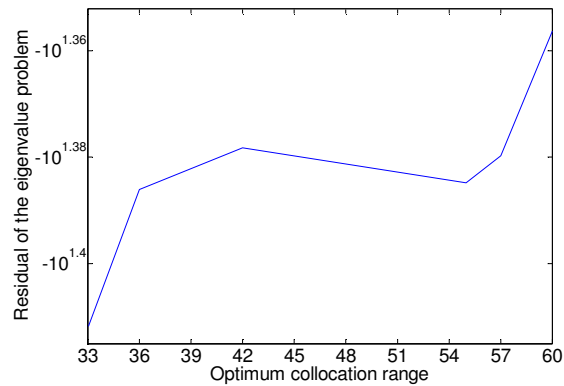

Fig.3 Residual along the optimum range for mode $m = 1$.

## 7 Conclusion

This paper reports two spectral approaches based on Chebyshev polynomials for numerically solving linear stability problems with applications to the swirling fluid system downstream the hydraulic Francis turbine. Numerical results showed that the use of these methods improved the computational time and the obtained critical values of the eigenparameter involved are accurate.

For bending modes the eigenvalue problem and its less simpler boundary conditions were translated into a linear system using a modified $L^2$ projection method and a collocation method, both based on shifted Chebyshev expansions. The numerical approximation of the unknown perturbation field was searched directly in the physical space. Both methods have been compared with the numerical results obtained in some other existing spatial investigations [29]. The collocation method proved to be more accurate, however the projection method was less expensive with respect to the numerical implementation costs, i.e. numerical results were obtained for a much smaller number of terms in the discretization. The results are very useful not only from their numerical values point of view, but also for their physical interpretation in fluid dynamics in flow control problems.